\documentclass[12pt,a4paper,amsmath,amssymb,graphicx,epsfig,subfigure]{article}
\usepackage{a4wide}
\usepackage{amsmath,amssymb, textcomp}
\usepackage{epsfig}
\usepackage{float}
\usepackage{comment}
\usepackage{textcomp}

\usepackage{tikz}
\usetikzlibrary{calc}

\tikzstyle{vertex} = [circle, draw=black, fill=black, inner sep=0.5mm, outer sep=0mm]
\tikzstyle{invvertex} = [circle, draw=white, fill=white, inner sep=0.0mm, outer sep=0mm]
\usetikzlibrary{shapes.geometric, arrows, positioning}
\usepackage{geometry}
\tikzstyle{title} = [rectangle, rounded corners, minimum width = 5cm, minimum height = 1cm, text centered, draw = black, fill=red!30]
\tikzstyle{mathematician} = [rectangle, rounded corners, minimum width = 3cm, minimum height = 1cm, text centered, draw = black, fill=red!10]

\tikzstyle{arrow} = [thick, ->, >=stealth]

\usepackage{float}
\usepackage{comment}
\usepackage{textcomp}
\usepackage{fancyhdr} 

\usepackage{hhline}

\newcommand{\BEA}{\begin{eqnarray}}
\newcommand{\EEA}{\end{eqnarray}}
\newcommand{\BEAs}{\begin{eqnarray*}}
\newcommand{\EEAs}{\end{eqnarray*}}

\begin{document}
\begin{center}
  \begin{large}\textbf{India and the Calculus of Trigonometric Functions}
  \end{large}
\end{center}
\vskip 0.5cm
\begin{center}
\textbf{Vijay A. Singh}\\
\textbf{Visiting Professor, Physics Dept., CEBS, Mumbai}\\
\textbf{Email: physics.sutra@gmail.com}\\ \vspace{0.1in}
\textbf{ABSTRACT}
\end{center}
We discuss the work of a brilliant line of Mathematicians who lived in central Kerala and starting with its founder Madhava (1350 CE) developed what can best be described as  Calculus and applied it to a class of trigonometric functions. We explain, with the example of the expansion of the inverse tan function, how they handled integration. Further, they took forward the work of the fifth century mathematician Aryabahata (499 CE), worked with differentials, and developed the expansion of the sine and the cosine functions. The work \textit{Yuktibhasa} (circa 1500 CE) which maybe described as the first textbook on Calculus, also describes in detail the evaluation of the area and volume of trigonometric functions as well as a variety of expansions for \textit{pi}. Our treatment is pedagogical and we present exercise problems and invite the enterprising student to try their hands at approximations and integration \textit{a la} the Madhava way.  

\section{Introduction}
Some six hundred years ago a  cluster of temple villages, on the banks
of the  Nila (now called  Bharathapuzha) river in central  Kerala, was
host to a brilliant line of mathematicians. Pre-eminent among them was
the founder  Madhava (1350 CE)  who pioneered  what came to  be called
Calculus. Little of  what the prescient Madhava wrote  has survived. A
lineage of  disciples not only kept  this flame of calculus  alive but
developed and wrote about it. It is this writing which is available to
us.  We mention  a few.  Parameshvara, who  was a  direct disciple  of
Madhava, wrote  profusely and spent  55 years examining the  night sky
and  documenting  eclipses.  He  along  with  his  two sons  Ravi  and
Damodara was a teacher to Nilakantha.  Among the many books Nilakantha
authored   the   \textit{Tantrasamgraha}   and   the   \textit{Bhasya}
(commentary)  on  a seminal  text  \textit{Aryabhatiya}  (499 CE)  are
notable.  His  student Jyeshthadeva is the  one we owe a  big debt to.
He wrote the \textit{Yuktibhasa}  which Divakaran (see References) has
designated  as  the ``first  text  book  on Calculus''.   The  lineage
continued till  the 1800s  and we  refer the reader  to Fig.1  and its
caption.

Madhava  and   his  students  developed  for   example  expansions  of
trigonometric  functions and  their inverses.   These expansions  were
developed a century or more later by European mathematicians using the
Calculus of  Newton and Leibniz. This  fact was noted and  reported by
Charles M.  Whish  \footnote{``On the Hindu Quadrature  of the Circle,
and the Infinite Series of the  Proportion of the Circumference to the
Diameter  Exhibited  in the  Four  Sastras,  Tantra Sanghraham,  Yucti
Bhasha,  Carana   Padhati,  and   Sadratnamala'',  Charles   M.  Whish
Transactions  of  the  Royal  Asiatic Society  of  Great  Britain  and
Ireland,  Vol.  \textbf{3},  pg.   509, (1834).   Whish knew  Shankara
Varman  (see Fig.1)  personally.}.   The Indian  written tradition  is
largely word  based.  Results  are mentioned  and the  derivations are
omitted.  The Aryabhatiya  (499 CE) with a little over  a 100 cryptic,
super-compressed verses of  dense mathematics is a  prime example.  In
the case of the Nila mathematicians however we are more fortunate.  We
can,  thanks particularly  to  the text  \textit{Yuktibhasa}, see  the
detailed reasoning although they are still word based. In what follows
we  shall provide  a  flavour  of the  methods  used  by these  Indian
mathematicians and  some of  their results. The  exercises in  the end
will help you get a more hands on understanding.

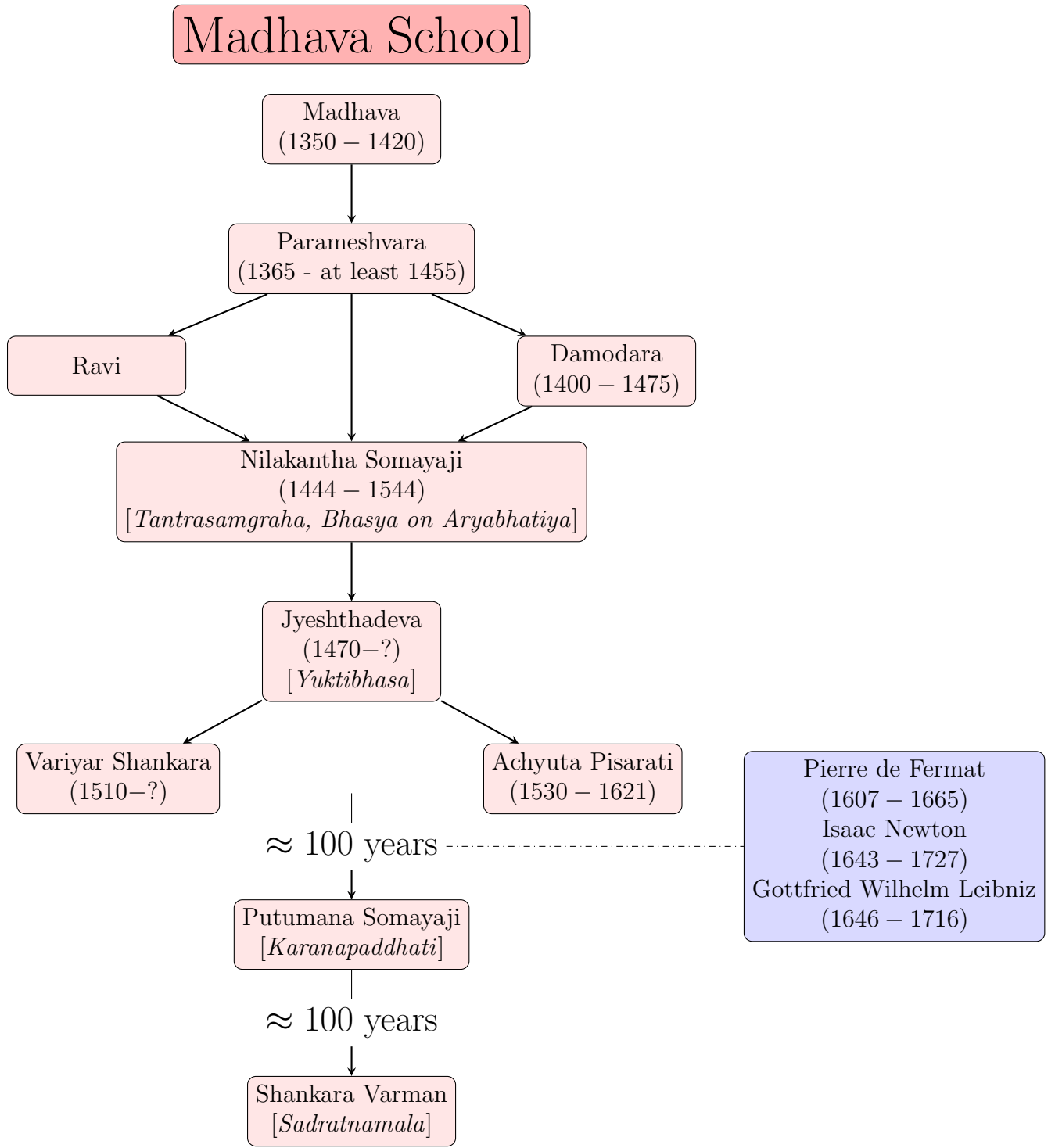
\begin{figure}
\begin{tikzpicture}[scale=.1]
\node (Nila School) [title] {\Huge Madhava School};
\node (Madhava) [mathematician, below = 0.5cm of Nila School, align=center] {Madhava \\($1350-1420$)};
\node (Parameshvara) [mathematician, below = 1cm of Madhava, align=center] {Parameshvara \\($1365$ - at least $1455$)};
\node (Damodara) [mathematician, below right = 1cm of Parameshvara, align=center] {Damodara \\($1400 - 1475$)};
\node (Ravi) [mathematician, below left = 1cm of Parameshvara, align=center] {Ravi};
\node (Nilakantha Somayaji) [mathematician, below = 2.5cm of Parameshvara, align=center] {Nilakantha Somayaji \\ ($1444 - 1544$)\\ $[$\textit{Tantrasamgraha, Bhasya on Aryabhatiya}$]$};
\node (Jyeshthadeva) [mathematician, below = 1 cm of Nilakantha Somayaji, align=center] {Jyeshthadeva \\ ($1470 - ?$)\\ $[$\textit{Yuktibhasa}$]$};
\node (Variyar Shankara) [mathematician, below left = 1cm of Jyeshthadeva, align=center] {Variyar Shankara \\ ($1510 - ?$)};
\node (Achyuta Pisarati) [mathematician, below right = 1cm of Jyeshthadeva, align=center] {Achyuta Pisarati \\ ($1530 - 1621$)};
\node (P) [circle, below = 1.5cm of Jyeshthadeva, align=center, fill=white, inner sep=0]{};
\node (Q) [draw=white, below = 0.5cm of P, fill=white, node distance=0.1cm, font=\Large]{$\approx 100$ years };
\node (Putumana Somayaji) [mathematician, below = 0.5 cm of Q, align=center] {Putumana Somayaji \\  $[$\textit{Karanapaddhati}$]$};
\node (R) [draw=white, below = 0.5cm of Putumana Somayaji, fill=white,  node distance=0.1cm, font=\Large]{$\approx 100$ years };
\node (Shankara Varman) [mathematician, below = 0.5cm of R, align=center] {Shankara Varman \\  $[$\textit{Sadratnamala}$]$};

\node (WestMath) [mathematician, right = 5.0 cm of Q, align=center, fill=blue!15] {Pierre de Fermat\\$(1607-1665)$\\Isaac Newton\\$(1643-1727)$\\Gottfried Wilhelm Leibniz\\$(1646-1716)$};

\draw [arrow] (Madhava) -- (Parameshvara);
\draw [arrow] (Parameshvara) -- (Damodara);
\draw [arrow] (Parameshvara) -- (Ravi);
\draw [arrow] (Parameshvara) -- (Nilakantha Somayaji);
\draw [arrow] (Damodara) -- (Nilakantha Somayaji);
\draw [arrow] (Ravi) -- (Nilakantha Somayaji);
\draw [arrow] (Nilakantha Somayaji) -- (Jyeshthadeva);
\draw [arrow] (Jyeshthadeva) -- (Variyar Shankara);
\draw [arrow] (Jyeshthadeva) -- (Achyuta Pisarati);
\draw (P) -- (Q);
\draw [arrow] (Q) -- (Putumana Somayaji);
\draw (Putumana Somayaji) -- (R);
\draw [arrow] (R) -- (Shankara Varman);
\draw [dashdotted] (Q) -- (WestMath);
\end{tikzpicture}
\caption{Members of the Madhava school.  Nilakantha's year of birth 1444 is firmly established. The other dates are tentative with some uncertainty ($\pm$ 5 years). Except for Ravi and Damodara who were sons of Parameshvara the other members are not direct descendants. The pioneering scholar of Indian mathematics Sarasvati Amma has designated this lineage as the \textbf{Aryabhata School}.}
\end{figure}


\section{\textit{Samskaram}: Recursive Refining}

As  a  methodology, recursion  has  been  used extensively  by  Indian
mathematicians.    It   would   be   best   to   explain   the   term
\textit{Samskaram} or recursive refining  with a simple example. There
is    another    name    for     it    -    \textit{Shudhikarna}    or
\textit{Shudikarti}. Consider  the evaluation  of $1/(x-d)$  given the
fact that we know $1/x$. We write
\begin{eqnarray*}
  \dfrac{1}{x-d} &=& \dfrac{1}{x} - \left[\dfrac{1}{x} - \dfrac{1}{x-d}\right] \\
                 &=&  \dfrac{1}{x} + \dfrac{d}{x}\left( \dfrac{1}{x-d}\right)
\end{eqnarray*}
On the r.h.s. of the second step we have the ``unknown term'' $1/(x-d)$. We replace it with the r.h.s. of step one,
\begin{eqnarray*}
  \dfrac{1}{x-d} &=&  \dfrac{1}{x} +  \dfrac{d}{x}\left[ \dfrac{1}{x} -  \left(\dfrac{1}{x} -  \dfrac{1}{x-d}\right)\right] \\
  &=& \dfrac{1}{x} + \dfrac{d}{x} \left[\dfrac{1}{x} +  \dfrac{d}{x} \left(\dfrac{1}{x-d}\right)\right]
\end{eqnarray*}  
and continuing recursively one more step
\BEAs
\dfrac{1}{x-d}  &=& \dfrac{1}{x} + \dfrac{d}{x} \left[\dfrac{1}{x} +  \dfrac{d}{x} \left(\dfrac{1}{x} - \left(\dfrac{1}{x} - \dfrac{1}{x-d}\right)\right)\right]\\
   &=& \dfrac{1}{x} + \dfrac{d}{x} \left[\dfrac{1}{x} +  \dfrac{d}{x} \left(\dfrac{1}{x} + \dfrac{d}{x}\left(\dfrac{1}{x-d}\right)\right)\right]
\EEAs
At this point we could drop the ``d'' on the extreme r.h.s to get
$$ \dfrac{1}{x-d} \approx \dfrac{1}{x} +  \dfrac{d}{x^2} +  \dfrac{d^2}{x^3} + \dfrac{d^3}{x^4}  $$
or ``refine'' our evaluation  further namely
\begin{eqnarray}
  \dfrac{1}{x-d} &=& \dfrac{1}{x} + \dfrac{1}{x} \sum_{n=1}^M \left(\dfrac{d}{x}\right)^n  \label{series1}
\end{eqnarray}
This method of iterative refining or recursive refining is called \textit{samskaram}. One of the earliest  usages of this method was for obtaining the square root of a number and is in the \textit{Bakshali} manuscript found near Peshawar and dated most probably 300 - 500 CE. The exercise at the end will give you a better feeling for this. We pause to note the following:
\begin{enumerate}
\item The method is not the same as the familiar Taylor expansion. In fact when used for the cosine series it yields
  $$cos(\theta +  \delta) \approx  cos(\theta) - \delta  sin(\theta) -
  \delta^2/2 \,\,cos(\theta) + \delta^3/8  \,\,sin(\theta) + .. $$ the
  last  term should  have  $\delta^3/6$ and  is  erroneously given  by
  \textit{samskaram}.  The Madhava school\footnote{To  avoid  confusion in  what
  follows we shall attribute all results to the ``Madhava school'' and
  only occasionally to the illustrious  ``Nila lineage'' or to a text
  such as  \textit{Yuktibhasa}. The members of  the lineage themselves
  from time to time use  the phrase ``\textit{Madhavoditam}'' (so said
  Madhava) and once in a while  invoke Aryabahata.} seemed to be aware
  of this and go on to derive the correct expansion (see Sec.\,IV).
\item For $x$ and $d$ positive and $d<x$ the series is convergent. The specific example cited by Nilakantha is $x=4$ and $d=1$.
\item The example above demonstrates a comfort level with infinite series. The l.h.s. is an unknown finite number and the r.h.s. is an infinite series which equals this number. This may not seem like an issue except when one views it from the perspective of Madhava. One is confronting for the first time an unknown irrational number $\pi$ and one claims to represent it with an infinite series - something that Madhava did. We shall see more of this later (Sec.\,III).
\end{enumerate}

  
\section{\textit{Samkalitam} - Integration}
\vskip -0.3cm  
\noindent \subsection{Introduction}
The term \textit{Samkalitam} means sum. In our case it is a special sum, one whose limit is an integral. The term will become clearer as we proceed in this section.

A crowning achievement of Madhava is the series representation of the angle $\theta$ in terms of $t=tan(\theta)$ for $\theta \leq \pi/4$.
  \begin{eqnarray}
 \theta &=& t- \dfrac{t^3}{3} + \dfrac{t^5}{5} - ... \,\,\,\,\, (t = tan(\theta)) \label{tanexp1}    
  \end{eqnarray}
  
We recognize this as the ``Gregory-Leibniz'' series and is one of the results which surprised Charles Whish  since it preceded European calculus by more than two centuries (see footnote Sec. I). We stress that it was meticulously derived and not just,  to use a cricketing terminology, a ``one-off lucky strike''. We shall derive this and thus get a taste of how the Nila lineage handled integration. Our demonstration proceeds in two stages. The first is the geometric part where Madhava, by intricate reasoning,  established the expression (Sec. III.2)
\begin{eqnarray}
  \delta \theta &=& \dfrac{\delta tan(\theta)}{1+ tan^2(\theta)}  \label{tandif1}
\end{eqnarray}
The next two sections deal with the \textit{Samkalitam} (integration) of the above expression (Sec III.2 and III.3). The entire procedure is a first in the history of Calculus.

\subsection{The Geometric Part}
Madhava obtained the relation between the angle $\theta$ and the tangent $t\,(=tan(\theta))$. Figure 2 depicts the quadrant of a unit circle. The circumscribing unit square is $OXWY$. We denote the angle $XOA_{n-1}$ as $\theta_{n-1}$ and the angle $XOA_{n}$ as $\theta_{n}$. The line $XW$ is divided into a large number $N$ of equal segments with $X = A_0, A_1, A_2, ..A_{n-1}, A_{n}, ... A_N = W$. The linear segment $A_{n-1} A_n$ = $1/N$ is small and equal to the increment $\delta t$ in the tangent. The corresponding increment in the arc of the unit circle is $PQ$ and equal to $\theta_n - \theta_{n-1}$. To repeat
\BEA    
\delta t &=& A_{n-1} A_n = 1/N             \label{exp1} \\
\delta \theta_n &=& \theta_n - \theta_{n-1}    \label{exp2} 
\EEA
Also from the right angle triangle $OXA_n$,
\BEA
  OA_n^2 &=& OX^2 + XA_n^2 \nonumber \\
         &=& 1 + (\dfrac{n}{N})^2   \label{exp3}
  \EEA

\begin{figure}
  \begin{center}
\begin{tikzpicture}[scale=0.8]

\draw [draw=white](10,10) grid (20,20);

\node [vertex, label=below:\Large \textbf{$O$}, name=O] at (10,10) {};
\node [vertex, label=below:\Large \textbf{$X$}, name=X] at (20,10) {};
\node [vertex, label=right:\Large \textbf{$W$}, name=W] at (20,20) {};
\node [vertex, label=left:\Large \textbf{$Y$}, name=Y] at (10,20) {};
\node [vertex, label=right:\Large \textbf{$A_{n}$}, name=A1] at (20, 19) {};
\node [vertex, label=right:\Large \textbf{$A_{n-1}$}, name=A2] at (20, 15) {};
\node [vertex, label=right:\Large \textbf{$P$}, name=P] at (18.94, 14.47) {};
\node [vertex, label=right:\Large \textbf{$Q$}, name=Q] at (17.43, 16.70) {};
\coordinate (r0) at ($(O)!(P)!(A1)$);
\coordinate (b0) at ($(O)!(A2)!(A1)$);
\node [vertex, label=left:\Large \textbf{$B$}, name=B] at (b0) {};
\node [vertex, label=left:\Large \textbf{$R$}, name=R] at (r0) {};
\node [invvertex, label=right:\Large \textbf{$\theta_{n}$}] at (12.5, 10.4) {};
\node [invvertex, label=right:\Large \textbf{$\theta_{n-1}$}] at (11, 10.3) {};

\draw (10,10) rectangle (20,20);
\draw  (10,10) -- (20, 19);
\draw  (10,10) -- (20, 15);

\draw (20,10) arc (0:90:10cm);
\draw (12.5,10) arc (0:42:2.5cm); 
\draw (11,10) arc (0:26.52:1cm); 
\draw (b0) -- (A2);
\draw (r0) -- (P);

\draw [rotate around=-48.1:(r0)](r0) rectangle ($(r0) + (0.2, 0.18)$);
\draw [rotate around=-48.1:(b0)](b0) rectangle ($(b0) + (0.2, 0.18)$);
\end{tikzpicture}    
\caption{A quadrant of the unit circle $OXY$ is circumscribed by a unit square $OXWZ$. The points $A_{n-1}$ and $A_n$ are very close to each other but the distance between them is magnified for ease of view. The arc PQ subtends an angle $\delta \theta_n = \theta_n - \theta_{n-1}$ at $O$. $PR$ and $A_{n-1} B$ are perpendiculars on $OA_n$ and these will be needed only for one of the Exercises.}
\end{center}
\end{figure}
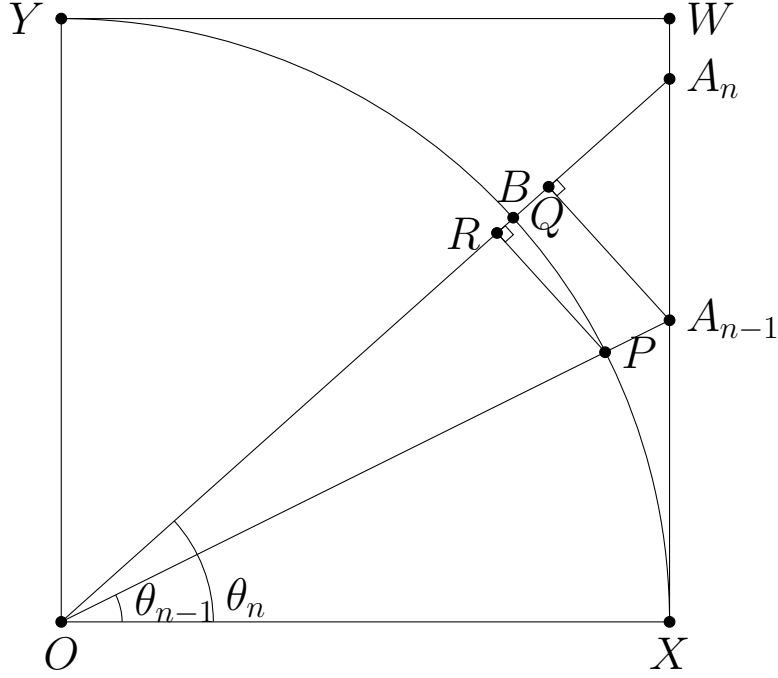

Through an elaborate series of arguments based on similar triangles and the smallness of $\delta \theta_n$ Madhava shows that
\BEA
sin(\delta \theta_n) &=& \dfrac{\delta t}{OA_n^2}       \nonumber \\
                    &=& \dfrac{1}{N (1 + (n/N)^2)}       \label{exp4}
\EEA
The \textit{Yuktibhasa} asks us to think of $N$ as very large; it uses the word
\textit{pararddham} or 10$^{17}$ and mentions that this is notional and to conceive of even larger numbers! The angle $\delta \theta_n$ and the $\delta t$ are, in its own words ``\textit{shunyaprayam}'' meaning of the nature of zero (and not ``\textit{shunyam}'' or zero). In other words this is the ``\textit{infinitesimal}'' of Calculus. To drive home the point the text also refers to it  as ``\textit{anuprayam}'' or atomic. Thus
\BEA
\delta\theta_n  &=& \dfrac{1}{N (1 + (n/N)^2)}       \label{exp5}
\EEA
Which is easily recognizable as Eq.(\ref{tandif1}). The next step is to integrate the expression Eq.(\ref{exp5}). We do it in two stages. 

\subsection{Samkalitam of the first few terms}
We can easily sum the l.h.s. of Eq.(\ref{exp5})
\BEA
\sum_{n=1}^N \delta \theta_n &=& (\theta_1 - \theta_0) + (\theta_2 - \theta_1) + ...\, (\theta_N - \theta_{N-1})  \nonumber \\
               &=& \pi/4    \label{exp51}
\EEA
We employ Eq.\,(\ref{series1}) with $x$ = 1 and $d$ = $-(n/N)^2$ to expand the r.h.s. of Eq.(\ref{exp5}) to obtain an infinite series
\BEA
\pi/4 & \approx & \dfrac{1}{N} \sum_{n=1}^N \left[1 - \dfrac{n^2}{N^2} + \dfrac{n^4}{N^4} - ..-\right]      \label{exp6} \\
& \approx & I_N(0) - I_N(2) + I_N(4) - I_N(6) + ... \\
\mbox{with} \nonumber \\
 I_N(k) & \approx & \dfrac{1}{N} \sum_{n=1}^N \left[\dfrac{n^k}{N^k}\right] \label{exp7} 
\EEA

The first three terms as displayed above are easily summed. 
\BEAs
I_N(0) = \dfrac{1}{N} \sum_{n=1}^N \left[1\right] = 1  \\
I_N(2) = \dfrac{1}{N} \sum_{n=1}^N \left[\dfrac{n^2}{N^2}\right]  = \dfrac{N (N+1) (2N+1)}{6 N^3} \\
I_N(4) = \dfrac{1}{N} \sum_{n=1}^N \left[\dfrac{n^4}{N^4}\right] =  \dfrac{N (N+1) (2N+1) (3N^2 + 3N - 1)}{30 N^5} \\
\EEAs
We next take the limit $N \rightarrow \infty$ and denote the limiting quantities by $J$. On inspection
\BEAs
 J_0 = \lim_{N\rightarrow \infty} I_N(0) = 1  \\
  J_2 = \lim_{N\rightarrow \infty} I_N(2) = 1/3  \\
 J_4 = \lim_{N\rightarrow \infty} I_N(4) = 1/5  \\
 \EEAs
  There are analytic expressions for $I_N(6)$ and $I_N(8)$ but higher orders would require a knowledge and manipulation of the Bernoulli numbers. Did Madhava and his disciples handle the higher orders, e.g. $J_{20}$ for example? The answer is they did so by induction and resorting to the asymptotic limit of large $N$. Their method is described in the \textbf{Appendix} and  represents yet another example of their mathematical acumen. Explicitly they obtain
 \BEA
 J_k &=& \frac{1}{k+1}          \label{exp8}
\EEA 
 Summarising, we have
 \BEA
 \dfrac{\pi}{4} &=& 1 - \dfrac{1}{3} +\dfrac{1}{5} - \dfrac{1}{7} + ... \nonumber\\
 &=& \sum_{k=0,2,4,...}^\infty (-1)^{k/2}\dfrac{1}{k+1}  \label{exp9}
 \EEA
 We close this section with a few pertinent remarks.
 \begin{enumerate}
 \item We can connect the above with the Calculus we have been taught. For example in Eq.\,(\ref{exp7}) above take
   \BEAs
   1/N \rightarrow dt ;\,\,\,\,\,\, \sum_1^N \rightarrow \int_0^1; \,\,\,\,\,\, (n/N)^k \rightarrow t^k    \\
   \EEAs
so   
\BEA   
   \dfrac{1}{N} \sum_{n=1}^N \left[\dfrac{n^k}{N^k}\right] \rightarrow \int_0^1 t^k dt
   \label{exp10}
   \EEA
 \item A key difference with the Calculus we are used to is that in the Madhava scheme we do \textbf{the summation first and then take the limit $N$ going to infinity}. In the former we take the limit first to get the differential $dt$ and then perform the integration. Since the derivatives of the powers of $t^k$ are known the fundamental theorem of Calculus\footnote{Namely the integral of the function $f$ over a fixed interval is the change in its anti-derivative F between the ends of the interval.} is invoked to mechanically write down the result $1/(k+1)$. In our case the integral is performed by first principles.  There are some advantages to the Madhava scheme as we shall see. One of them, namely the interchange of summations, is mentioned above. The fundamental theorem follows trivially in our summation procedure. An example is  Eq.(\ref{exp51}) where the summation is replaced by the end points ($\pi/4$ -$0$). We will point this out with another example in the next section. 
 \item We can take $\delta t = tan(\theta)/N$ instead of  $\delta t = 1/N$ resulting in the general series
   \BEA
   \theta &=& 1 - \dfrac{tan^3(\theta)}{3} +\dfrac{tan^5(\theta)}{5} - \dfrac{tan^7(\theta)}{7} + ...     \label{exp(11)}
   \EEA
   We can legitimately use the term ``function'' here. It describes the dependence of a real quantity $\theta$ on another real variable $\tan{\theta}$. In other words we have a functional expression for every value of the variable.   
   
 \item Prior to Madhava we had an approximate value for $\pi$, namely 22/7 or  as given by Aryabhata, namely 62832/20000 = 3.1416. Aryabhata is clear that this  value is \textit{assana} i.e. proximate,\footnote{This word is to be distinguished from \textit{sthula} or approximate. For example when theorists use the value of the mass of the electron to be $9.1\times10^{-31}$ kg  it is \textit{sthula}. When experimentalists carefully quote the value 9.10938 $\times$ 10$^{-31}$ kg it is \textit{assana} meaning that it can be refined with more careful experimentation.}. meaning that it is close to but not quite $\pi$. The point to appreciate is that instead of another proximate value, Madhava has an \textbf{exact} infinite series for $\pi$. This suggests the irrationality of $\pi$ but Madhava does not clearly say so.
The series in Eq.(\ref{exp9}) is a slowly converging one and a number of methods were proposed to develop rapidly convergent series The world record up to the 1800s was held by Shankara Varman (see Fig. 1) of the Madhava school with $\pi$ up to 18 decimal places. 
\end{enumerate}

 
\section{Differentials; Sine and Cosine Expansions}
 In this Section we initially dwell on the work of Aryabhata and then see how it led the Madhava school to the develop the notion  of the differential and the expansions of the Sine and Cosine series.
 \subsection{The Aryabhata Connection}
 One can discern a continuity in Indian mathematics, howsoever tenuous,  from pre-Vedic times ($<$ 1000 BCE) up and until 1800s. A striking example is the influence of the  \textit{Aryabhatiya} (499 CE) on the Madhava school (1350 CE).

 The \textit{Aryabhatiya} has some 121 verses out of  which 33 verses belong to the mathematics section (\textit{Ganitapada}). Aryabhata defines, for the first time in the history of mathematics, the sine function. It is the half-chord $AP$ of the unit circle in Fig. 3.
\BEAs
sin(\theta) &=& \dfrac{AP}{OA}  \\
            &=& AP   \,\,\,\,\,\,\,\, (OA = 1)
\EEAs
The circle maybe large or small; correspondingly $AP$ and $OA$ maybe large of small, but the l.h.s. is a function of $\theta$ and is \textbf{invariant}. With this, Aryabhata endowed circle geometry with metrical properties. This alone may qualify him as the founder of trigonometry. But he did more.

\begin{figure}
\begin{center}  
\begin{tikzpicture}[scale=.5]

\draw [draw=white](10,10) grid (20,20);
\draw (15,15) circle (5cm);
\draw (15,15) -- (15, 20);
\draw (15,15) -- (20, 15);
\draw (15, 15) -- (18, 19);
\draw [style=dashed](18, 19) -- (18, 15);
\draw (16, 15) arc (0:24.5:2.5cm);
\node [vertex, label=below:\Large \textbf{$O$}, name=O] at (15,15) {};
\node [vertex, label=above:\Large \textbf{$Y$}, name=Y] at (15,20) {};
\node [vertex, label=right:\Large \textbf{$X$}, name=X] at (20,15) {};
\node [vertex, label=above:\Large \textbf{$A$}, name=A] at (18,19) {};
\node [vertex, label=below:\Large \textbf{$P$}, name=P] at (18,15) {};
\node [invvertex, label=right:\Large \textbf{$\theta$}, name=theta] at (15.9,15.5) {};
\end{tikzpicture}
\caption{The quadrant of the unit circle. The half chord $AP$ is $sin(\theta)$ as defined by Aryabhata. See text for comments.}
\end{center}
\end{figure}
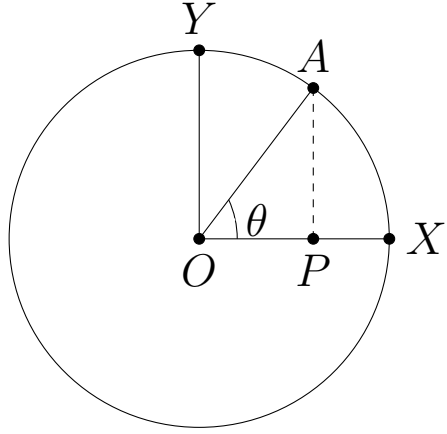

He knew that the difference in the sines of two angles $\phi + \delta \phi$ and $\phi - \delta \phi$ is proportional to the cosine of the mean angle $\phi$, 
\BEA
sin(\phi + \delta \phi) - sin(\phi - \delta \phi) &=& 2 sin(\delta \phi) cos(\phi) \label{sindiff1}
\EEA  
and the difference in the cones of two angles $\phi + \delta \phi$ and $\phi - \delta \phi$ is proportional to the sine of the mean angle $\phi$.
\BEA
cos(\phi + \delta\phi) - cos(\phi - \delta \phi) &=& - 2 sin(\delta \phi) sin(\phi) \label{cosdiff1}
\EEA
Further he states the second sine difference. According to his commentators this is done by (once again)  ingenuous arguments based on similar triangles. We shall not go over it here since our concern is different here. But Aryabhata did more.

He then obtained the values of the sines at fixed angles between  0 and $\pi /2$ thus  generating the sine table for $\pi/48$ = 3.75 degrees, 7.5 degrees ... 90 degrees. This table has been used by Indian astronomers (and astrologers) in some form or another since 499 CE up to the present. We shall see this in the next section. To give you an idea every time you use the calculator, or look up a table, to search for the values of the  sines and trigonometric functions, you may silently thank Aryabhata for showing the way. But Aryabhata did more. 

\subsection{Finite Difference Calculus}
Let us take $\delta \phi = \epsilon$  where $\epsilon$ is small but not infinitesimal (not ``\textit{shunyaprayam}''). We take $\phi = n \epsilon$ where $n$ is a positive integer from 1 to $N$.  To fix our ideas $\epsilon = \pi$/48 = 3.75$^0$ = 3438'. Employing the sine and cosine difference formulae from the previous section we define differences
\BEA
\delta s_{n} &=& s_{n+1} - s_{n} = 2 s_{1/2} c_{n+1/2} \label{finite1} \\
 \delta c_n &=& c_{n} - c_{n-1} = -2 s_{1/2} s_{n+1/2} \label{finite2} 
 \EEA
 where the symbol $s_n$ stands for $sin(n \epsilon)$, $c_n$ stands for $cos(n \epsilon)$ and $s_{1/2}$ for $sin{\epsilon/2}$. Note that these are essentially the same as Eqs. (\ref{sindiff1}) and (\ref{cosdiff1}), e.g. $sin(\phi + \delta \phi) - sin(\phi) = 2 sin(\delta \phi/2) cos(\phi +\delta \phi/2)$). The above is a pair of coupled equations and it was Aryabhata's insight to take the second difference, namely
\BEA
\delta^2 s_n &=& \delta s_{n} - \delta s_{n-1}  = 2 s_{1/2} (c_{n+1/2}-c_{n-1/2})  \nonumber  \\
& =& - 4 s_{1/2}^2 s_n \mbox{\,\,\,\,\,\,\, using Eq.\,(\ref{finite2}) }  \label{finite3}
 \EEA 
 Thus the second difference of the sines is proportional to the sine itself. Similarly the second difference in the cosine is also proportional to the cosine.
 \BEA
\delta^2 c_n &=& \delta c_{n} - \delta c_{n-1} = 2 s_{1/2} (s_{n+1/2}-s_{n-1/2})  \nonumber  \\
& =& - 4 s_{1/2}^2 c_n \mbox{\,\,\,\,\,\,\, using Eq.\,(\ref{finite1}) }  \label{finite4}
 \EEA 

 The next step is to represent the r.h.s in terms of a recursion. We observe $s_n$ on the r.h.s. of Eq.(\ref{finite3}) may be written as $s_n = s_n - s_{n-1} + s_{n-1} - s_{n-2} + s_{n-2} - ...$ = $\delta s_{n-1} + \delta s_{n-2} + ...$  Thus
\BEA 
 \delta s_{n} - \delta s_{n-1} &=& - 4 s_{1/2}^2 \sum_{m=1} ^{n-1} \delta s_m \label{recur1}
 \EEA

 A few remarks are in order here.
 \begin{enumerate}
 \item   The above work is that of Aryabhata and he takes $\epsilon$ = $\pi$/48. He also approximates $s_1 = sin(\epsilon) \approx \epsilon$. Using it he derived the sines of angles from 3.75 $^0$ to 90$^0$ in 24 equi-spaced steps. We will suggest a simple problem along these lines in the Exercise.
\item  Of more relevance is the fact that the above is the algorithm we would currently use to obtain derivatives numerically. We know that sine(37$^0$) is close to 0.6 and sine(30$^0$) is 0.5. Thus the difference in angle is 7$^0$ which in radians is 0.122. Thus derivative of sine of the median angle 33.5$^0$ is from Eq.\,(\ref{finite1}) is
 \BEAs
 \delta sin(\phi) /\delta \phi  &=& (0.6 - 0.5)/0.122 = .82
 \EEAs
 Looking up the sine table or the calculator yields cos(33.5) = 0.83. Similarly Eq.\,(\ref{finite3}) yields the second derivative namely
 \BEAs
   \delta^2 sin(\phi) / \delta^2 \phi  \approx - sin(\phi)  
 \EEAs
 The above are now called central difference approximations to the derivative and the second derivative. Aryabhata does not mention the term  finite difference calculus (let alone calculus). But similar methods  are now  used to numerically solve our differential equations. That includes Newton's II Law and the famous Schrodinger equation of quantum mechanics both of which are  second order differential equations. 
\end{enumerate} 
 
\subsection{The Trigonometric Series}
Aryabhata took $\delta \phi = \pi$/48 (= 3.75$^0$ = 3438'). But this is not sacrosanct and he states that $\delta \phi$ could be ``\textit{chindyat yateshtani}'' meaning as small as you like. Or as large. Brahmagupta (600 CE) for example took it to be $\pi$/12 to generate the sine series using Aryabhata's formulae from the previous section. To his credit he also had robust methods to evaluate the sine for intermediate values. The Madhava school looked in the other direction. They took it to be very small.

We rewrite Eqs.\,(\ref{sindiff1}) and (\ref{cosdiff1}) with $\delta \phi$ replaced by
 $\delta \phi/2$.
\BEA
\delta sin(\phi) &=& sin(\phi + \delta \phi/2) - sin(\phi - \delta \phi/2)  \nonumber \\
                 &=& 2 sin(\delta \phi/2) cos(\phi) \label{sindiff2} \\
\delta cos(\phi) &=& cos(\phi + \delta\phi) - cos(\phi - \delta \phi) \nonumber \\
                 &=& - 2 sin(\delta \phi/2) sin(\phi) \label{cosdiff2}
\EEA
Like in the previous section we take $\delta \phi/2 = \phi/N$ where once again $N$ is unimaginably large, larger than \textit{paraddham} or 10$^{17}$! Then $\delta \phi/2$ is \textit{sunyaprayam} or ``of the nature of zero'' i.e., an infinitesimal. We can then take replace 2 $sin(\delta \phi/2)$ by $\delta \phi$ to yield
\BEA
 \delta sin(\phi) &=& \delta \phi \,cos(\phi) \label{sindiff3} \\ 
 \delta cos(\phi) &=& - \delta \phi \,sin(\phi) \label{cosdiff3} 
 \EEA
 The Madhava school leaves the above equations in the differential form. They do not explicitly write the derivative. But nothing prevents us from doing so.
 \BEA
 \lim_{\delta \phi \rightarrow 0} \dfrac{ \delta sin(\phi)}{\delta \phi} &=& \dfrac{d sin(\phi)}{d \phi} \nonumber \\
 &=& cos(\phi) \label{sindiff4}
 \EEA
 Similarly
 \BEA
 \dfrac{d cos(\phi)}{d \phi} &=& - sin(\phi)   \label{cosdiff4}
\EEA
Taking inspiration from the previous section  we  may take one more derivative to obtain
e.g. $\delta^2 sin \phi$ and  $\delta^2 cos \phi$
\BEA 
\dfrac{d^2 sin(\phi)}{d \phi^2} = - sin(\phi)  \label{sindiff5} \\
\dfrac{d^2 cos(\phi)}{d \phi^2} = - cos(\phi)  \label{cosdiff5} 
\EEA
Madhava and his disciples worked in the  discrete domain leaving the limiting procedure for the end. We shall illustrate how they proceeded using the modern notation so familiar to us. 
Integrating once from zero to $\theta$
\BEAs
 \left. \dfrac{d sin(\phi)}{d \phi}\right|_{\theta} - \left. \dfrac{d sin(\phi)}{d \phi}\right|_{0} &=& - \int_0^{\theta} sin(\phi) d \phi \,\,\,\,\,\,\,\mbox{Or}    \\
\left.  \dfrac{d sin(\phi)}{d \phi}\right|_{\theta} &=& 1 -  \int_0^{\theta} sin(\phi) d \phi
\EEAs  
Note that we have used the fundamental theorem of calculus. We repeat this procedure once more to obtain
\BEA
 sin(\theta) &=& \theta - \int_{0}^{\theta} d\phi \int_0^{\phi} sin(\xi) d \xi  \label{sams}
\EEA
The above expression can be readily subjected to recursive refining (``\textit{samskaram}''). To start with we have
$$ \sin(\theta)_1 = \theta$$
Next we replace the sine in the second term on the r.h.s. of Eq.(\ref{sams}) by $sin(\xi) = \xi$. Hence for the second recursion
\BEAs
sin(\theta)_2 &=& \theta -  \int_{0}^{\theta} d\phi \int_0^{\phi} \xi d \xi \\
&=&  \theta - \int_0^{\theta} \phi^2/2\, d\phi    \\
 &=&  \theta - \theta^3/3!    
\EEAs
It is easy to see the trend. Next we replace the $sin(\xi)$ in Eq.\,(\ref{sams}) by $\xi -\xi^3/3!$. This yields
\BEAs
sin(\theta)_3 &=&    \theta - \theta^3/3! + \theta^5/5!
\EEAs
The entire sine series is thus obtained
\BEA
sin(\theta) &=& \sum_{k=1,3,5...} (-1)^{(k-1)/2}\dfrac{\theta^k}{k!}   \label{sams1}
  \EEA
  One may similarly obtain the cosine series
  \BEA
  cos(\theta) &=& \sum_{k=0,2,4...} (-1)^{k/2} \dfrac{\theta^k}{k!}    \label{sams2}
  \EEA
  We have carried out the expansion using the method of ``\textit{Samskaram}'' but interpolated with the more familiar integration and appeal to the fundamental theorem of calculus. In what one may describe as a \textit{tour de force}, the book \textit{Yuktibhasa} carries out the entire exercise in the discrete formalism describing it in minute detail in the  Malayalam language\footnote{Most manuscripts by Indian mathematicians are in Sanskrit and in verse form respecting the norms of grammar and meter of Sanskrit}. We have forgone this. The Appendix which pertains to the previous section will give one a feeling for how this discrete formalism works.  

  
  \section{Discussion}
  Among the other accomplishments of the Madhava school we mention two. They derived the formulae for the area and the volume of the sphere by methods of calculus. Secondly they realized that the $\pi$ series (Eq.\ref{exp9}) has slow convergence. So they reformulated the series in multiple ways. A proof of their ingenuity is the calculation of the value of $\pi$ up to 18 decimal places by Sankara Varman (1800s, see Fig.1). At that time this was a world record. We shall not describe these accomplishments here.

  From time to time one one hears of the Calculus accomplishments of Indian mathematicians predating Madhava. While describing the motion of celestial bodies, Bhaskara II (1100CE)  has used terms like \textit{tatkalika} (at that instant). He also mentions the stationarity property of elliptical orbits at the apogee and perigee. It is a creditable example of theoretical insight based on observation. A close reading reveals that he is still thinking in terms of small and finite differences in angles (not time). Similar claims have been made on behalf of Nilakantha (see Fig.\,1) with his refined astronomical model. This is as per the Aryabhata framework described in Sec.IV.2. Quantities are small, but there is no infinitesimal (``\textit{shunyaprayam}'') and the careful treatment it requires.  The area and volume of the sphere are also mentioned but the methods by which they are arrived at are unclear. There is a parallel in European mathematical history. Archimedes arrived at the value of $\pi$ by the method of exhaustion long ago. Both Descartes and Fermat had discussed ``derivatives'' in terms of the slope. Fermat even mentions points in the function $f(x)$ where a small change in $x$ has ``almost no effect'' on f(x). One can describe the work of these illustrious mathematicians as pre-Calculus or proto-Calculus at best.

  Some shortcomings are apparent in the work of the Madhava school. They were sensitive to the convergence of infinite series. But they did not prove the convergence, absolute or conditional. Neither did Newton. An explicit derivative notation or its interpretation in terms of slope is not present. In  a sense the Madhava school's treatment of differentials is closer to Leibniz than Newton. One can sense their reluctance to divide ``zero by zero''. There is no treatment of conic sections (hyperbola for instance). Further, how would one accomplish the integration of say $\sqrt{x}$ in \textit{Samlkalitam}? Or of exponential and logarithmic functions? The answer to these questions is that perhaps it is possible. It is for students, particularly Indian students,  to push forward  the Madhava program and to explore its advantages and limitations. All told it is a beautiful approach to Calculus. The Madhava school was in decline by the time serious European science and mathematics came to India and perhaps did not have a chance to address these concerns.

To sum up, the Madhava school had a consistent formalism using methods that can be identified as methods of Calculus and they applied it successfully to a class of functions, namely, trigonometric functions. The approach is refreshingly different from the European.  To quote the Fields Medallist David Mumford, ``It seems fair to me to compare [Madhava] with Newton and Leibniz''.\footnote{Notices of the American Mathematical Society, Vol. \textbf{57}, page 385 (2010).}   

\paragraph{Acknowledgement:} The author would place on record the many useful discussions he had with Prof. P. P. Divakaran. 


\appendix
\section{ Derivation of $J_k = 1/(k+1)$ (Eq.\,(\ref{exp8}))}
  We define a related quantity $S_N(k)$
  \BEA
  S_N(k) = N^{k+1} I_N(k) = \sum_{n=1}^N n^k                \label{a1}
  \EEA
  The proof proceeds in two parts. We first obtain a recursion relation for $S_N(k)$. Next we take the asymptotic limit of large $N$ to obtain an explicit expression for  $S_N(k)$.

  To obtain the recursion relation we first note that \BEA N S_N(k-1)
  - S_N  (k) &=&  \sum_{n=1}^{N} (N  -  n) n^{k-1}  \label{a2} \\  &=&
  \sum_{n=1}^{N-1} \sum_{j=1}^n j^{k-1} \label{a3}  \EEA Thus \BEA S_N
  (k) &=&  N S_N(k-1) - \sum_n^{N-1}  S_n (k-1) \label{a4} \EEA  It is
  not  easy to  see how  the  single summation  of Eq.\,(\ref{a2})  is
  reordered to the double summation of Eq.\,(\ref{a3}). One way to see
  this  is to  take  some  concrete values  say  $N=5$  and $k=3$  and
  convince oneself  thereby. Divakaran's  book (see  References) takes
  the continuum limit  of this and uses integration by  parts to prove
  it. He  also points  out that  this is  a special  case of  the Abel
  re-summation formula  and it is remarkable  that \textit{Yuktibhasa}
  had   discovered   and   deployed    it.    Using   the   definition
  (Eq.\,(\ref{a1})) we can  arrive at the recursion  relation given by
  Eq.\,(\ref{a4}) from Eq.\,(\ref{a3}).

The next step is to solve for $S_N (k)$ in the large $N$ limit. We begin by noting from Eq.\,(\ref{a1}) that $S_N(0) = N$. Hence  Eq.\,(\ref{a4}) yields    
\BEA
S_N (1)   &=&   N^2 - \sum_n^{N-1} S_n (0) \nonumber  \\
           &=&   N^2 - (1 + 2 + 3 + ... (N-1))  \nonumber \\
           &=& N^2 - N(N-1)/2      \nonumber \\
 &=&      N^2/2      \label{a5}  
\EEA
where we take the large N limit in the last step.  Eq.\,(\ref{a5}) suggests that
\BEA
S_N (k) = N^{k+1}/(k+1)      \label{a6}
\EEA
We then use mathematical induction. We shall establish that $S_N(k+1) = N^{k+2}/(k+2)$ using the recursion formula  Eq.\,(\ref{a4}). Note
\BEA
S_N (k+1)   &=&   N S_N(k) - \sum_n^{N-1} S_n  (k)  \nonumber \\
&=& N^{k+2}/(k+1)  - \sum_n^{N-1} n^{k+1}/(k+1)    \nonumber \\
&=&  \dfrac{N^{k+2}}{k+1}  - \dfrac{ S_{N-1}(k+1)}{k+1}  \label{a7}
\EEA
where the last step follows from the definition of $S_N(k)$ in  Eq.\,(\ref{a1}).
Since N is large we take  $S_{N-1}(k+1) \approx  S_N(k+1)$. This then establishes our required relation
\BEA
 S_N (k+1)   &=&  \dfrac{N^{k+2}}{k+2}     \label{a8}
 \EEA
We now note that we have proved this asymptotically. Hence  we replace $I_N (k)$ in Eq.\,(\ref{a1}) by the $N \rightarrow \infty$  limiting  expression $J_k$, 
\BEA
J_{k+1} & = & \dfrac{1}{k+2}    \label{a9}
\EEA
and similarly $J_k = 1/(k+1)$.  We note that the derivation holds for all positive integers $k$ and not just for even integers that were required in Sec. III.

\vskip 1.0cm
\begin{center}
  \textbf{EXERCISES}
\end{center}
\begin{enumerate}
\item \textit{Samskaram} for square root of a positive number $n$. Take a perfect square $m^2$ less than but closest to $n$.
  Define a correction $n = m^2 +r$.  We can rewrite this a
  $$ n = (m+r/2m)^2 - (r/2m)^2. $$
  For the first iteration drop the $(r/2m)^2$ term on the r.h.s. so
  $$ \sqrt{n}_1 = m +r/2m.$$
  Continue this iteration and show that
  $$ \sqrt{n}_2 = m + r/2m - (r/2m)^2 \dfrac{1}{2(m + r/2m)}.  $$
  Numerically compute for $n$ = 95 and $m$ = 9.
\item An alternate \textit{Samskaram} for the square root of the positive number $n$. Take
  $$ \sqrt{n}_1 = \dfrac{1}{2} \left(m + \dfrac{n}{m}\right) .$$
  Obtain the next iteration. Once again numerically compute for the $n$ = 95 and $m$ = 9. Which of the two methods yields a closer value for $\sqrt{5}$?
\item \textit{Samskaram} method for the cosine interpolation formula of Madhava. Consider the two trigonometric identities:
  \begin{eqnarray}
    cos(\theta +\delta) &=& cos(\theta) - 2 sin(\delta/2) sin(\theta +\delta/2)   \label{arya1} \\
    sin(\theta +\delta) &=& sin(\theta) + 2 sin(\delta/2) cos(\theta +\delta/2)   \label{arya2}
  \end{eqnarray}
  We point out that the Nila lineage acknowledge that these identities were first mentioned in the Ganita section of Aryabhatiya (499 CE). Madhava then approximates $sin(\delta/2)$ by $\delta/2$ for small $\delta$ to write
\begin{eqnarray}  
   cos(\theta +\delta) &=& cos(\theta) - \delta sin(\theta +\delta/2)   \label{arya3} 
\end{eqnarray}  
  and next drops the $\delta/2$ in the sine term on the r.h.s. to obtain the first step of the recursion
  $$ cos(\theta +\delta)_1 = cos(\theta) - \delta sin(\theta)   $$
  If on the other hand we had retained Eq.(\ref{arya3}) and employed the exact sine formula Eq.(\ref{arya2}) we would refine the recursion. Show that continuing one will obtain
  $$ cos(\theta +\delta)_3 = cos(\theta) - \delta sin(\theta) -\frac{\delta^2}{2} \,cos(\theta) + \dfrac{\delta^3}{8}\, sin(\theta)   $$
  This does not lead to the correct expansion for the cosine. The Nila mathematicians were well aware of the limitations of the above expansion. The point of this exercise is to point out that the method of \textit{Samskaram} is often accompanied by other approximations (in this case the $sin(\delta/2) \approx \delta/2$), whose reliability must be gauged independent of the iterative process.  
\item The Indian mathematical tradition made judicious use of similar triangles. We can get a taste of this by establishing  the important relation between the  angle and its tangent (Eq.(\ref{exp4})). First prove that the triangles $OXA_n$ and $A_{n-1}A_nB$ are similar and further that triangles $OPR$ and $OA_{n-1}B$ are similar.  Thus establish
  $$ PR = \dfrac{A_nA_{n-1}}{OA_n\,OA_{n-1}}   $$
  Note that $PR = sin(\delta \theta_n)$ and $OA_n \approx OA_{n-1}$ while $A_nA_{n-1} = \delta t$.

\item We can generate the sine table as per Aryabhata's suggestion but not exactly using the same value for $\epsilon$ he used. We choose $\epsilon = \pi/80 \approx 0.0393$ which is the same as 4.5$^0$. We take $sin(\epsilon) \approx \epsilon$ and  $sin(2 \epsilon) \approx 2 \epsilon)$. If you have a simple calculator generate all values of sine from 2.25 to 18 degrees in equal steps using Eq.\,(\ref{recur1}). Alternatively if you have a programmable calculator or a computer generate all values of sine from 2.25 to 90 degrees. Compare with the results your calculator will otherwise yield.      

\item Using the method of Section IV.3 generate the cosine expansion to arrive at
  Eq\,.(\ref{sams2}).

\item Take the continuum limit of the summation formulae (Eq.\,(\ref{a2}) and Eq.\,(\ref{a3})). To do this see Eq.\,(\ref{exp10}). Prove that the two summations are indeed equal. [Note: This proof is due to Divakaran (see References).]     
  
\end{enumerate}


\begin{center}
  \textbf{FEW REFERENCES}
\end{center}
These references are in English. \\
\textbf{Primary}
\begin{enumerate}
\item  ``Aryabhatiya'', Walter E. Clark, University of Chicago Press (1930).
\item ``Aryabhatiya'', Kripashanker Shukla and K V Sarma, Vols. I and II, Indian National Science Academy Publication (1976). The Section IV.1 and IV.2 have been shaped by the critical discussions of certain verses from the \textit{Ganitapada} by these authors. 
\item ``Ganita-Yukti-Bhasa'' (of Jyeshthadeva), K. V. Sharma (with commentary by K. Ramasubramanian, M. D. Srinivas, and M. S. Sriram), Hindustan Book Agency, New Delhi (2008).
\end{enumerate}

\textbf{Secondary:}
\begin{enumerate}
\item ``The Mathematics of India'', P. P. Divakaran, Hindustan Book Agency (2018). The book has shaped this article in ways covert and overt. Highly recommended.   
\item ``The History of Hindu Mathematics's'', Bibhutibhusan Datta and Avadhesh Narayan Singh, Vols. I and II, Asia Publishing House, Delhi (1935 and 1938). A pioneering book on Indian Mathematics written in Pre-Independence India.
\item ``Geometry in Ancient and Medieval India'', Sarasvati Amma, Motilal Banarsidas, Delhi (1999). It has detailed discussions worth looking at. It also uncovers the element of continuity in the Indian mathematical traditions from ancient times to the pre-British era.     
\end{enumerate}
\vskip 1.0cm

\textit{Prof. Vijay A. Singh  has been faculty at IIT Kanpur (1984-2014) and HBCSE, Tata Institute for Fundamental Research (2005-2015) where he was the National Coordinator of  both the Science Olympiads and the National Initiative on Undergraduate Science for a decade. He is a Fellow, National Academy of Sciences, India and was President, the Indian Association of Physics Teachers (2019-21). He is curruently a Visiting Professor CEBS, Mumbai. (emailid: physics.sutra@gmail.com)}
\end{document}